\documentclass[12pt]{article}
\usepackage{latexsym}
\usepackage{amssymb,amsmath,graphicx}
\usepackage{amsfonts}
\usepackage{epsfig}
\input epsf
\def\RR{{\mathbb{R}}}
\def\CC{{\mathbb{C}}}

\def\p{\partial}
\def\z{\zeta}

\def\Rea{\mathrm{Re}\,}
\def\Ima{\mathrm{Im}\,}

\begin{document}
\title{Quasi-exceptional domains}
\author{Alexandre Eremenko\thanks{Supported by NSF grant DMS-1361836.}$\;$
and Erik Lundberg}
\maketitle
\begin{abstract}
Exceptional domains are domains on which there exists a positive harmonic
function, zero on the boundary and such that the normal derivative
on the boundary is constant.
Recent results classify exceptional domains as belonging to either 
a certain one-parameter family of simply periodic domains or one of its scaling limits.

We introduce quasi-exceptional domains by allowing the boundary
values to be
different constants on each boundary component.
This relaxed definition retains the interesting property
of being an \emph{arclength quadrature domain},
and also preserves the connection to the hollow vortex
problem in fluid dynamics.
We give a partial classification of such domains in terms of certain
Abelian differentials.
We also provide a new two-parameter family of periodic quasi-exceptional domains.
These examples generalize the hollow vortex array found by Baker, Saffman, and Sheffield (1976).
A degeneration of regions of this family provide doubly-connected examples.

2010 AMS Subject Class: 35R35, 76B47, 30C20, 31A05

\end{abstract}

\noindent
{\bf 1. Introduction}
\vspace{.1in}

A domain $D\in\RR^n$ is called exceptional if there is
a positive function $u$ (called a \emph{roof function}) harmonic in $D$,
zero on the boundary, and
\begin{equation}\label{1}
\frac{\partial}{\partial n}u(z)=1,\quad z\in\partial D ,
\end{equation}
where the differentiation is along the normal pointing inwards into $D$,
and it is assumed that the boundary is smooth.
Evident examples are exteriors of balls and half-spaces.
For $n>2$ the only other known examples are cylinders
whose base is an exceptional domain in $\RR^2$.
If the smoothness assumption on the boundary is dropped, then
there are also certain cones in higher dimensions
and pathological ``non-Smirnov'' examples in the plane \cite{HL}.

The problem of description of all exceptional domains in the plane was stated
in \cite{HHP} and settled in \cite{HL} under a topological assumption which was removed in \cite{MT}
using an unexpected correspondence to minimal surfaces. 
The first non-trivial example was given in \cite{HHP}.
This example appeared in another context related to fluid dynamics in \cite{L}.
A second non-trivial example was noticed in \cite{HL} and \cite{MT}.
This example had also appeared previously in studies of fluid dynamics \cite{1976} (see also \cite{Crowdy}).

Let us introduce \emph{quasi-exceptional} domains, by relaxing the definition to allow the Dirichlet condition to be a different constant on each boundary component.
Thus, a domain $D\in\RR^n$ is called quasi-exceptional 
if there is a positive harmonic function $u$ in $D$,
which is constant on each boundary component (but not necessarily the same constant) and
the Neumann condition (\ref{1}) holds (with the same constant on each boundary component).
We will continue to call $u$ a \emph{roof function}.
Again, we assume that each component of the boundary is smooth.

We summarize several interesting aspects of exceptional domains.
These statements all hold true for quasi-exceptional domains except the last one:

\begin{itemize}
 \item Fluid dynamics: As noted above, the two non-trivial examples first appeared in fluid dynamics \cite{L, 1976}.
In general, one can interpret exceptional domains in terms of a \emph{hollow vortex} problem.
 The level lines of $u$ can be interpreted as stream lines
of a two-dimensional stationary flow of ideal fluid, and condition (\ref{1}) expresses the fact
that the pressure is constant on the boundary. Such conditions may exist
if the components of the complement of $D$ are air bubbles in
the surrounding liquid.  Notice that the rotation of the
fluid around all bubbles corresponding to exceptional domains is
in the same direction. This reflects our condition that $\p_n u>0$.
\item Quadrature domains \cite{Gu}: Exceptional domains provide examples of arclength null-quadrature domains, that is, domains for which integration over $\partial D$
of every analytic function in the Smirnov class $E^1(D)$ vanishes.
\item Differentials on Riemann surfaces: By way of the connection to quadrature domains, the study \cite{Gu}
indicates a connection to half-order differentials.
We make use of Abelian differentials in Section 4 below.

\item The Schwarz function of a curve: In \cite{HL}, it was noticed that the function $u(z)$ satisfies $$ \partial_z u(z) = \sqrt{-S'(z)}, $$
where $S(z)$ is the Schwarz function of $\p \Omega$ and $\partial_z = \frac{1}{2} (\p_x - i \p_y )$ is the Cauchy-Riemann operator.
\item Minimal surfaces: The recent work \cite{MT} established a nontrivial correspondence between exceptional domains
and a special type of complete embedded minimal surfaces called a ``minimal bigraphs''.
This correspondence does not extend to quasi-exceptional domains.

\end{itemize}

The classification results for exceptional domains show that they are quite restricted;
all examples can be conformally mapped from a disk by elementary functions.

\medskip

\noindent {\bf Problem 1:} Classify quasi-exceptional domains.

\medskip

We begin to address this problem below, give a partial classification
of periodic and finitely-connected exceptional domains,
and provide new periodic and doubly-connected examples
described in terms of elliptic functions.
First, we explain the relation to arclength null-quadrature domains.

\medskip

\noindent {\bf 2. Arclength null-quadrature domains}

\medskip

A bounded domain $D \subset \CC$ is a \emph{quadrature domain} if it admits a formula expressing the area integral of 
every function $f$ analytic and integrable in $D$ as a finite sum of weighted point evaluations of the function and its derivatives.
i.e. 
\begin{equation}\label{eq:QF1}
 \int_{D} {g(z) dA(z)} = \sum_{m=1}^{N} \sum_{k=0}^{n_m}{a_{m,k}g^{(k)}(z_m)},
\end{equation}
where $z_m$ are distinct points in $D$ and $a_{m,k}$ are constants independent of $g$.

A (necessarily unbounded) domain $D \subset \CC$ is called a \emph{null-quadrature domain} (NQD) if the area integral of 
every function $g$ analytic and integrable in $D$ vanishes:
\begin{equation}\label{eq:QF2}
 \int_{D} {g(z) dA(z)} = 0.
\end{equation}
M. Sakai \cite{Sakai1981} completely classified NQDs in the plane.

Following \cite{HL} we will refer to a domain $D \subset \CC$ as an \emph{arclength null-quadrature domain} (ALNQD) if the integral over $\partial D$ of 
every function $g$ in the Smirnov class $E^1(D)$ vanishes 
(in the case $\infty$ is an isolated point on $\p D$, 
we take the restricted class of functions $g(z) \in E^1(D)$ vanishing at infinity):
\begin{equation}
 \int_{\partial D} {g(z) ds(z)} = 0.
\end{equation}

The Smirnov class $E^1(D)$ is not the same as the Hardy space $H^1(D)$.
Namely, a function $g$ analytic in $D$ is said to belong to $E^1(D)$ if there exists a sequence of cycles $\gamma_k$
homologous to zero, rectifiable, and converging to the boundary $\p D$ 
(in the sense that $\gamma_k$ eventually surrounds each compact sub-domain of $D$),
such that:
$$\sup_{\gamma_k} \int_{\gamma_k} |g(z)| |dz| \leq \infty.$$

One may also define quadrature domains in higher dimensions using a test class of harmonic functions, 
but we will restrict ourselves to the case of $n=2$ dimensions.

Inspired by the successful classification of NQDs \cite{Sakai1981}, 
the problem of classifying ALNQDs was suggested in \cite{HL}.
We pose this problem again while stressing that it does not reduce to the classification of exceptional domains
(whereas it might reduce to classification of \emph{quasi-exceptional} domains).

\medskip

\noindent {\bf Problem 2:} Classify ALNQDs.

\medskip

The following proposition shows that quasi-exceptional domains are ALNQDs.
Thus, the new examples (described in the last section) of quasi-exceptional domains also provide new ALNQDs.
Problem 2 is closely related to Problem 1, and
if the converse of the proposition is true then the two problems are equivalent.

\medskip

\noindent
{\bf Proposition 1.} {\em Suppose that $D$ is a quasi-exceptional domain.
Then $D$ is an ALNQD.}

\medskip

\noindent {\em Proof.}
Consider the complex analytic function $F(z) = u_x - iu_y$, where $u$ is the roof function.
We will need the following Claim which is proved in the next section (see Lemma 2).

\smallskip

\noindent {\bf Claim.} {\em  The roof function $u$ of $D$ satisfies $\nabla u(z) = O(1)$ in $D$. Thus, $F(z)$ is bounded. }
\vspace{.1in}

Suppose that $g$ is in the Smirnov space $E^1(D)$, and let $\phi$ be a conformal map to $D$ from a circular $\Omega$.
Using the fact that $ds = i F(z) dz$, we have
\begin{equation}\label{eq:null}
\int_{\p D} g(z) ds = \int_{\p D} i g(z) F(z) dz = i \int_{\p \Omega} g(\phi(\z)) F(\phi(\z)) \phi'(\z) d\z.
\end{equation}
Now $g(z) \in E^1(D)$ implies \cite{Duren} that
$$g(\phi(\z)) \phi'(\z) \in E^1(\Omega),$$ 
unless $\infty$ is an isolated boundary point
(in which case $g(\phi(\z)) \phi'(\z)$ need not be analytic at $\phi^{-1}(\infty)$).
In the case $\infty$ is not an isolated boundary point, since $F(\phi(\z))$ is bounded (by the Claim), we have 
$$g(\phi(\z)) F(\phi(\z)) \phi'(\z) \in E^1(\Omega),$$
and therefore, by Cauchy's theorem, (\ref{eq:null}) vanishes.
In the case $\infty$ is an isolated boundary point of $D$, $g(\phi(\z)) \phi'(\z)$ need not be analytic at $\z = \phi^{-1}(\infty)$.
In fact it can have a pole up to order two.
However, $F(z) = O(|z|^{-1})$ vanishes to order one at infinity by B\^ocher's Theorem (cf. \cite[Thm. 3.1]{HL}).
Thus, if $g(z)$ also vanishes at infinity then $g(\phi(\z)) \phi'(\z) F(\phi(\z)) \in E^1(\Omega)$ and the previous argument follows.
This shows that if $z \in D$ then $D$ is an ALNQD for the restricted class of functions $g(z) \in E^1(D)$ which vanish at infinity.
\vspace{.1in}

\noindent {\bf 3. A potential theoretic restriction on the roof function}

\vspace{.1in}

We restrict ourselves to the case $n=2$, and
assume that the order of connectivity of $D$
is finite, or that the roof function $u$ is periodic,
and the fundamental region for $D$ has finite connectivity. 

Recall that a Martin function is a positive harmonic function
$M$ in a domain $\Omega$
with the property that for any positive harmonic function
$v$ in $\Omega$, the condition $v\leq M$
implies that $v=cM$, where $c>0$ is a constant.
(Often, Martin functions are called minimal harmonic functions - cf. \cite{Heins1950}.)
Martin functions on finitely connected domains are simply Poisson kernels evaluated at points of the Martin boundary, 
the boundary under Caratheodory compactification (prime ends) of the domain (see \cite{Brelot71}).

Any domain $D$ of finite connectivity in $\CC$ is conformally equivalent
to a circular $\Omega$, whose boundary components are circles and points.
For a circular, a Martin function $M$ can be of two types:

a) There is a component of $\partial\Omega$ which is a single point
$z_0$, and $M$ is proportional to the Green function of $\Omega\cup\{ z_0\}$
with the singularity at $z_0$.

b) There is a point $z_0\in\partial\Omega$ which is not a component
of $\partial\Omega$, and $M$ has zero boundary values at all
points of $\partial\Omega\backslash\{ z_0\}.$ The local behavior
in this case is like $-\Ima(1/z)$ in the upper half-plane near $0$.

Let $D$ be an exceptional domain,
and $u$ a harmonic function with the property
(\ref{1}).
The following result was proved for exceptional domains by the current first author, but communicated in \cite[Thm. 4.2]{HL}.
Here we repeat the proof with minor adjustments.
\vspace{.1in}

\medskip

\noindent
{\bf Lemma 2.} {\em  The roof function $u$ of a quasi-exceptional domain satisfies $\nabla u(z) = O(1)$ in $D$.
Moreover, $u$ is a sum of a bounded harmonic function and at most two Martin functions.}
\vspace{.1in}

{\em Proof.}
We follow the second part of the proof from \cite[Thm. 4.2]{HL}.

Let $R>0$ and consider an auxiliary function
$$w_R=\frac{|\nabla u|}{u+R},$$ where $R>0$ is a parameter.
A direct computation shows that
\begin{equation}\label{dwa}
\Delta\log w_R = w_R^2,
\end{equation} 
and $w_R(z)=1/(c_k + R) \leq 1/R$ for $z\in\partial D$,
where $c_k \geq 0$ are the constants taken in the Dirichlet condition.
We claim that 
\begin{equation}\label{tri}
w_R(z)\leq 2/R,\quad z\in D,
\end{equation}
from which the result follows by letting $R \rightarrow \infty$ which gives $|\nabla u| \leq 2$ in $D$.

Suppose, contrary to (\ref{tri}), that $w_R(z_0)>2/R,$ for some $z_0\in D$.
Let 
$$v(z)=\frac{2R}{R^2-|z-z_0|^2},\quad z\in B(z_0,R)=\{ z:|z-z_0|< R\}.$$
Obviously, $v(z)\geq 2/R$.  A computation reveals that
$\Delta\log v=v^2.$ 
Let
$$K=\{ z\in D \cap B(z_0,R): w_R(z)>v(z)\}.$$
We have $z_0\in K$, since $v(z_0)=2/R$.
Let $K_0$ be the component of $K$, containing $z_0$. Then we have
$w_R(z)=v(z)$ on $\partial K_0$, since $w_R(z)< v(z)$ on $\partial D \cap B(z_0,R)$
while $v(z)=+\infty$ on $\partial B(z_0,R)$.
On the other hand,
$$ \Delta ( \log w_R-\log v ) \geq w_R^2-v^2>0 \quad\mbox{in}\quad K_0.$$
So the subharmonic function $\log w_R -\log v$ is positive in $K_0$ and vanishes on the boundary,
a contradiction.

This proves that $\nabla u = O(1)$.
In order to see the second statement, we note that $\nabla u = O(1)$
implies that $u(z) = O(|z|)$ has order one.
The result then follows by first solving the Dirichlet problem (with a bounded function)
having the same boundary values as $u$; subtracting this function, one may then apply \cite[Theorem II]{Kjellberg}.
\vspace{.1in}

\noindent {\bf 4. Partial classification
in terms of Abelian differentials}
\vspace{.1in}

Let $D$ be a QE domain of one of the following types:

Type I. $D$ is finitely connected, or

Type II. $D/\Gamma$ is finitely connected,
where $\Gamma$ is the group of transformations $z\mapsto z+n\omega$,
and $u(z+\omega)=u(z)$ for some $\omega\in\CC\backslash\{0\}$.
We call this the periodic case. (As above, $u$ is the roof function.)

In this section we give a classification of QE domains of these two types
in terms of Abelian differentials of a compact Riemann surface with
an anti-conformal involution.

If $D$ is of type I, and $\infty$ is an isolated boundary point, then
$D'=D\cup\{\infty\}$ is conformally equivalent to some bounded
circular domain $\Omega$, and we suppose that $p \in \Omega$ corresponds
to $\infty$. If $\infty$ is not isolated, we put $D'=D$, and
$\Omega$ is a bounded circular domain conformally equivalent to $D'$. In any case,
we have a conformal map $\phi:\Omega\to D'$.

If $D$ is of type II, let $G=D/\Gamma$. The Riemann surface $G$ is a finitely
connected domain on the cylinder $\CC/\Gamma$; this cylinder
is conformally
equivalent to the punctured plane; $G$ must have one or two punctures
of $\CC/\Gamma$ as isolated boundary points, and we denote by $G'$ the union
of $G$ with these isolated boundary points. Then $G'$ is conformally equivalent
to a bounded circular domain of finite connectivity $\Omega$ and
we have a multi-valued conformal map $\phi:\Omega\to G'$. Let $a$ and $b$ 
denote the one or two points in $\Omega$ that correspond to the added punctures of $G'$.

In all cases $a$ and $b$ are simple poles of $\phi$.

We pull back $u$ on $\Omega:$ set $v=u\circ \phi$. As $u$ is periodic,
$v$ is a single-valued positive harmonic function on
$\Omega\backslash\{ a,b\}$.

Consider the differential on $\Omega$
$$dv=v_zdz=(1/2)(v_x-iv_y)(dx+idy)=g(z)dz.$$
This is well defined on $\Omega$:
$g$ is a single-valued meromorphic function
in $\Omega$ with simple poles exactly at $a$ and $b$ (if any
of these points is present in $\Omega$).

Next, we extend $v$ as a multi-valued function to a compact Riemann surface $S$.
Let $\Omega'$ be the mirror image of $\Omega$;
we glue it to $\Omega$ in the standard way (along each circular boundary component)
and obtain a compact Riemann surface $S$.
We denote by $\sigma:z\mapsto z^*$ the anti-conformal involution which fixes the boundary components of $\Omega$. 
The Riemann surface $S$ is
of genus $g$, and the involution $\sigma$ has fixed set corresponding to 
$\partial \Omega$, which consists of $n=g+1$ ovals. 
Such involutions are called involutions of maximal type.

Each branch of $v$ is constant on each boundary component,
so it extends through this boundary component by reflection to the
double $S$ of $\Omega$.
The extensions of various branches of $v$ through different
boundary components do not match: they differ by additive constants.
On the other hand, the differential $dv$ is well defined on the double.
Namely, 
\begin{equation}\label{diff-symmetry}
(dv)^* = -dv,
\end{equation}
where $*$ is the action of involution on differentials.
Thus we have a meromorphic differential $dv$ on $S$.

Choose a basis of $1$-homology in $S$ so that the $A$-loops 
are simple closed curves in $\Omega$, each homotopic to one boundary
component of $\Omega$, and the $B$ loops are dual to the $A$-loops.
For Type I, all periods over $A$-loops are purely imaginary, because
$$v=\Rea\int dv$$
is single-valued. For Types II these periods are imaginary except those which
correspond to  simple loops around one pole, $a$ or $b$.

Now we discuss $\phi$, or better the differential $d\phi=\phi'(z)dz$.
We have, from the condition that our domain is quasi-exceptional:
$$2|dv|=|d\phi|.$$
The ratio of two differentials is a function. So we have a meromorphic function $B$
on $\Omega$ such that
\begin{equation}\label{eq:B}
2Bdv=d\phi.
\end{equation}
This function has absolute value $1$ on $\partial\Omega$. Therefore,
it extends to $S$ by symmetry. Its poles belong to $\Omega$ and must match
the zeros of $dv$, because $d\phi$ is zero-free (indeed, $\phi$ is univalent).
In fact, $B$ is a meromorphic function on $S$.
To justify this claim when $dv$ has a singularity on $\p \Omega$,
we observe that this singularity is removable for $B$
as follows from the next lemma.
\vspace{.1in}

\noindent
{\bf Lemma 3.}
{\em Consider the equation
$$\phi'=Bh,$$
where $h$ is meromorphic in a neighborhood $V$ of $0$,
$B$ is holomorphic and zero-free in $V \backslash \{0\}$, 
$|B(z)|=1$ for $z \in V \cap \RR \setminus \{0\}$,
and $\phi$ is univalent in $\{ z\in V:\Ima z>0\}$.
Then the singularity of $B$ at $0$ is removable.}
\vspace{.1in}

Before proving the lemma, we note that in order to apply it in our setting
we compose $B$ with a linear fractional transformation that sends $V$ to a neighborhood of the singularity we wish to remove
such that the real line is mapped to the circular boundary component with $0$ sent to the singularity.

\vspace{.1in}

\noindent {\it Proof.}
Proving this by contradiction, assume that $0$ is an essential singularity
of $B$. By symmetry we have $B(\overline{z})=1/\overline{B(z)}$.
Then by Phragm\'en--Lindel\"of theorem, there exists
a sequence $z_k\to 0$ such that
\begin{equation}\label{eq:seq1}
\liminf_{k\to\infty}|z_k|\log|B(z_k)|>0.
\end{equation}
By symmetry, there exists a sequence $z_k^\prime\to 0$ such that
\begin{equation}\label{eq:seq2}
\liminf_{k\to\infty}|z_k^\prime|\log|B(z_k^\prime)|<0.
\end{equation}
Without loss of generality, we may assume that one of these sequences
$z_k$ or $z_k^\prime$ is in the upper half-plane.

Distortion theorems for univalent functions imply that
\begin{equation}\label{dist}
c(\Ima z)^3\leq |\phi'(z)|\leq C(\Ima z)^{-3},
\end{equation}
In addition we have
\begin{equation}\label{rea}
c|z|^m \leq |\phi'(z)| \leq C|z|^{-m}, \quad z \in V \cap \RR,
\end{equation}
for some $m>0$.
These two inequalities imply via Carleman's ``loglog'' principle \cite{C,R}
that
$$c|z|^m\leq|\phi'(z)|\leq C|z|^m,\quad z\in V,\quad\Ima z>0.$$
This contradicts either (\ref{eq:seq1}) or (\ref{eq:seq2}), depending on
which sequence $z_k$ or $z_k^\prime$ lies in the upper half-plane.
\vspace{.1in}

We can thus restate the problem of finding QE domains
(under the restrictions we impose) as follows:
\vspace{.1in}

{\em
Find a triple $(S,d\omega,B)$, 
where $S$ is a compact Riemann surface with an involution
of ``maximal type'' (the complement of the fixed set of the involution consists of
two regions homeomorphic to plane regions), $d\omega$ is a meromorphic
differential that enjoys the symmetry property (\ref{diff-symmetry}),
and $B$ is the function which has the symmetry
property
$$B^*(z):=\overline{B(z^*)}=1/B(z),$$
and has poles at the zeros of $d\omega$ on one half of $S$,
that is in $\Omega$.
There is an additional condition that 
$$\phi=2\int Bd\omega$$
is globally univalent and single-valued in case I, and
single-valued except the residues in case II.}
\vspace{.1in}

In order to check the
condition on the global univalence of $\phi$,
it is sufficient to verify that periods of
$d\omega/B$ are zero on the boundary curves,
and that these boundary curves are mapped by $\phi$ injectively.

A general conclusion is the following.
\vspace{.1in}

 \noindent
{\bf Proposition 4.} {\em The boundary of a quasi-exceptional domain
of type I or II is parametrized by an Abelian integral.}
\vspace{.1in}

Next we provide a partial classification of quasi-exceptional domains in terms of the data stated in the above formulation.

\medskip

\noindent
{\bf Theorem 5.} {\em The differential $dv$ has either two or four poles in $S$ counting multiplicity.
Moreover, if $dv$ has two poles in $S$ then $D$ is either a disk or a half-plane.}

\vspace{.1in}

\smallskip

\noindent {\em Remark.} If $B\neq {\mathrm{const}}$,
then $1/B$ an Ahlfors function of $\Omega$.

\medskip

\noindent {\em Proof.}
The differential $dv$ has simple poles at $p$, $a$, and $b$ (when present)
and at their images $\sigma p$, $\sigma a$, and $\sigma b$. 
In addition it may have double poles on $\partial\Omega$.
The total number of poles in $\overline{\Omega}$ is at most two
by Lemma~2.
Thus on $S$, the differential $dv$ has two or four poles,
counting multiplicity.

Notice that $v$ is constant on each boundary component,
so the gradient is perpendicular to the boundary $\partial\Omega$,
so the total rotation
of this gradient, as we describe the boundary is the same as the total
rotation of the tangent vector to the boundary. This is equal to 
$2\pi(2-n)$ because $C_1$ is traversed counterclockwise and the rest
clockwise, as parts of the boundary of $\Omega$.
So $v_z$ which is conjugate to the gradient, rotates $n-2$ times.

From this we can conclude how many zeros $dv$ has in $\Omega$.
The number $N$ of zeros of $dv$ in $\Omega$ satisfies
\begin{equation}
\label{count}
n-2 = N - (\mbox{the number of poles in}\; \Omega),
\end{equation}
where a double pole on $\partial\Omega$ is counted as a single pole in
$\Omega$.

Suppose $dv$ has exactly $2$ poles, counting multiplicity.  This can occur in one of three ways:

(1) $dv$ has a simple pole at $p$ in $\Omega$.

(2) $dv$ has one double pole at $z_0 \in \p \Omega$.

(3) $dv$ has a simple pole at $a$ in $\Omega$ (and $b$ does not exist).

If Case (1) holds,
then $\infty$ is an isolated point on $\p D$, and by Proposition 1, $D$ is an arclength quadrature domain with quadrature point at $\infty$.
It now follows from \cite[Remark 6.1]{Gu} that $D$ is a disk.

In Case (2), we will show that $B$ is constant.
First note that $d \phi$ has a double pole at $z_0$,
so $B$ does not have a zero or a pole at $z_0$.
Since $\phi$ is a conformal map, it follows from (\ref{eq:B}) that $B$ has no zeros and $N$ poles in $\Omega$ (located at the zeros of $dv$).
Assume for the sake of contradiction that $B$ is not constant.
By Lemma 3, $B$ is meromorphic in $S$, and by Lemma 2, $1/|B|$ is bounded by a constant in $\Omega$.
Since $|B|=1$ on $\p \Omega$, $B$ thus maps $\Omega$ to the exterior of the unit disk and maps each of the $n$ components of $\p \Omega$ to the unit circle.
This implies that $B$ has at least $n$ poles in $\Omega$.
Combined with (\ref{count}), this gives the contradiction $N = n-1 \geq n$.
We conclude that $B$ is constant which implies that the gradient of the roof function is constant.
Thus, the roof function is linear, and $D$ is a halfplane.

In Case (3), the behavior of $\phi$ at point $a$ is logarithmic, so $d \phi$ has a simple pole at $a$ and $B$ does not have a zero or a pole at $a$.
Arguing as before, we conclude that $B$ is constant and that $D$ is a halfplane.

\medskip

\noindent
{\bf Corollary 6.} {\em The only quasi-exceptional domain $D$ with compact boundary is the exterior of a disk,
and the only quasi-exceptional domain for which $\infty$ is a limit point of only one component of $\partial D$ is the halfplane.}

\vspace{.1in}

If $D$ is a quasi-exceptional domain that is not a disk or halfplane,
then $dv$ has four poles and more precisely, we have one of two possibilities:

\smallskip

$D$ is of {\em type I}: $dv$ has two double poles
on $\partial\Omega$.
This implies that the boundary $\partial D$ consists of two 
simple curves tending to $\infty$ in both directions,
and $n-1$ bounded components.
The unbounded components are the $\phi$-images of two arcs of
of one boundary circle of $\Omega$ which contains both singularities
of $\phi$ and $v$.

\smallskip

$D$ is of {\em type II}: $dv$ has two simple poles in $\Omega$.
In this case $D$ must be periodic, all components of $\partial D$ are
compact and there are $n$ such components per period.
\smallskip

Note that the possibility that 
$dv$ has one simple pole in $\Omega$ and
one double pole on $\partial\Omega$ is excluded by Lemma 2:
it is easy to see that in
this case the number of Martin functions in the decomposition of $u$ would be
infinite.

We have thus described possible topologies of the QE domains
satisfying the assumptions stated in the beginning of this section.

In the next section we construct
the examples of types I and II with $S$ of genus $1$.
We conjecture that there exist QE domains of types I and II
based on $S$ having any genus.
\vspace{.1in}

\noindent
{\bf 5. New examples}
\vspace{.1in}

Description of our examples
requires elliptic functions
(all known exceptional domains can be parametrized by elementary functions).

\vspace{.1in}

\noindent
{\bf Example of type I.}
\vspace{.1in}

Let $G$ be the rectangle with vertices $(0,2\omega_1,2\omega_1+\omega_3,
\omega_3)$, where $\omega_1=2\omega$, $\omega>0$,
and $\omega_3=\omega^\prime$,
where $\omega'\in i\RR,$ $\omega'/i>\omega$.
Let $G'$ be the reflection of $G$ in the real line.
The union of $G,G'$ and the interval $(0,2\omega_1)$ make a fundamental
domain of the lattice $\Lambda$ generated by $2\omega_1,2\omega_3$.

Let us consider the $\omega_1$-periodic positive harmonic
function $h$ in $G$ which is
zero on the horizontal segments of the boundary $\partial G$,
except for one singularity per period, at $0$, 
where it behaves in the following way:
$$h(z)\sim-\Ima(1/z),\quad z\to 0.$$
Note that the existence of $h$ is clear as it can be expressed
(through conformal mapping) in terms of the Poisson kernel of a ring domain.

Function $h$ has two critical points in $G$, at $w_1$ and $w_2$ with
$\Rea w_1=\omega_1/2$ and $\Rea w_2=3\omega_1/2$, while the imaginary
parts of $w_1,w_2$ are equal. Let us choose real constants $c_1$ and $c_2$
such that
$v=2(h+c_1y)+c_2$ is a positive harmonic function with
critical points $\omega_1/2+\omega_3/2$ and $3\omega_1/2+\omega_3/2$.
The existence of such constants $c_1$ and $c_2$ is evident by continuity.

The $z$-derivative $\p_z v = (v_x - i v_y)/2$ is 
an elliptic function with periods $\omega_1,2\omega_3$, and thus also
elliptic with periods $\Lambda$. 
Asymptotics near $0$ show that $\p_z v \sim -i/z^2$, and
as this function has only
one pole per period, (with respect to the parallelogram $\omega_1,2\omega_3)$,
we have $\p_z v = -i \wp + i c_0$, where $\wp$ is the Weierstrass
function corresponding to the lattice $(\omega_1,2\omega_3)$. 

Zeros of $\p_z v$ in $G\cup G'$ are
$\omega_1/2+\omega_3/2,3\omega_1/2+\omega_3/2$ and complex conjugates
in $G'$. 

Let $B$ be an elliptic function with periods $2\omega_1,2\omega_3$
having simple poles at $\omega_1/2+\omega_3/2,3\omega_1/2+\omega_3/2$,
and zeros at complex conjugate points. Such function exists by Abel's
theorem: the sum of zeros minus the sum of poles equals $-2\omega_3$.
This function is unique up to a constant factor.
By symmetry, $B(\overline{z})=c/\overline{B(z)},$
so on the real line $|B(x)|^2=c$ and we can choose the constant factor
in the definition of $B$ so that $c=1$. Thus 
\begin{equation}\label{111}
|B(x)|=1, \quad x\in\RR.
\end{equation}
Then we have $B(x+\omega_3)\overline{B(x-\omega_3)}=1$,
but by periodicity we also have $B(x+\omega_3)=B(x-\omega_3)$,
thus $|B(x+\omega_3)|=1$. So 
\begin{equation}\label{222}
|B(z)|=1\quad\mbox{on the horizontal segments of}\quad\partial G.
\end{equation}
Now we consider the function
$$F=\frac{\partial v}{\partial z}B= (-i\wp+ic_0) B.$$
This function $F$ is holomorphic and zero-free in $G$ (the zeros of 
$\partial v/\partial z$ in $G$ are exactly canceled by the poles of $B$).
Let us show that
\begin{equation}\label{34}
\int_0^{2\omega_1}F(x+iy)dx=0, \quad y \in(0,\omega_3).
\end{equation}
This property follows from the fact that $B(z)$ and $B(z+\omega_1)$ have
the same poles but the residues at these poles are of the opposite signs,
because $B$ has only two poles in the period parallelogram.
Thus
\begin{equation}\label{sim2}
B(z+\omega_1)=-B(z).
\end{equation}
Property (\ref{sim2}) and $\omega_1$-periodicity of $\wp$ imply (\ref{34}).

We conclude that the primitive $f=\int F$ is locally univalent.
Assuming for the moment that it is univalent,
it maps $G$ onto some region in the plane, and  we have
$$|f'|=|F|=\left|\frac{\partial v}{\partial z}\right||B|.$$
Define $u$ by composing $v$ with $f^{-1}$, so $u(f(z)) = v(z)$. 
Then $u$ is positive and harmonic in $f(G)$.
Taking into account (\ref{222}), we conclude that $u$ satisfies (\ref{1})
$f(G)$ is a quasi-exceptional domain.
Note that, in accordance with the previous results in \cite{MT}, 
$f(G)$ is not an exceptional domain since the piecewise constant Dirichlet data is not the same constant on each boundary component.

\begin{figure}[ht]
\centering
    \includegraphics[scale=.2]{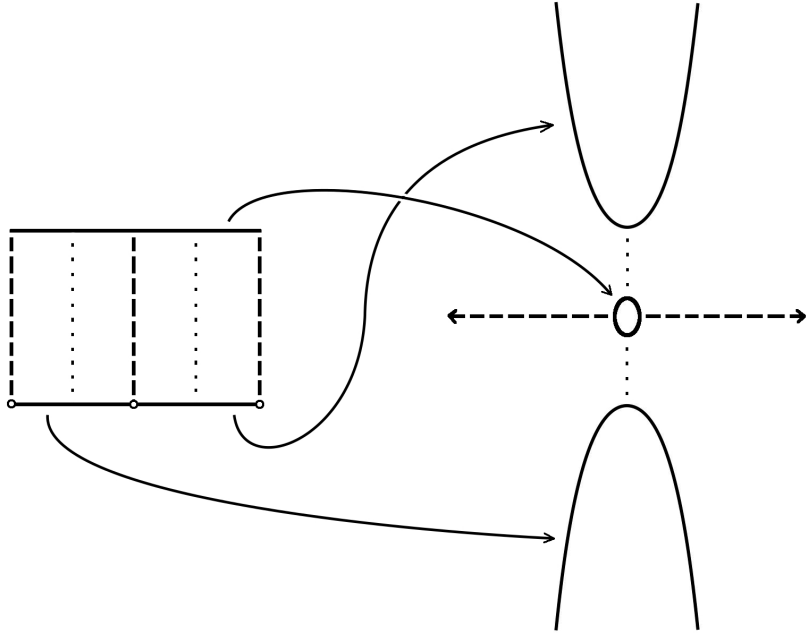}
    \caption{A doubly-connected quasi-exceptional domain of type I mapped from the rectangle $G$.}
    \label{fig:ExampleC}
\end{figure}
In order to show that $f$ is in fact univalent, it is enough to show that it is one-to-one on the the horizontal sides of $G$ (since $f$ is locally univalent).
To this end, we make the following claims:

\smallskip

\noindent Claim 1: $\Rea f$ is increasing along the segment $[\omega',\omega' + 2\omega]$ and decreasing along the segment $[\omega'+2\omega,\omega' + 4\omega]$.

\smallskip

\noindent Claim 2: $\Ima f < \Ima f(\omega')$ on $[\omega',\omega' + 2\omega]$ and $\Ima f > \Ima f(\omega')$ on $[\omega'+2\omega,\omega' + 4\omega]$.

\smallskip

\noindent Claim 3: $\Ima f$ achieves its minimum and maximum values along $[\omega',\omega' + 4\omega]$ at $\omega' + \omega$ and $\omega' + 3 \omega$ respectively.

\smallskip

\noindent Claim 4: $\Rea f$ is increasing on the segment $[0,2 \omega]$, and $\Rea f$ is decreasing along $[2 \omega, 4 \omega]$.

\smallskip

\noindent Claim 5: $\Ima f$ attains its maximum along $[0,2 \omega]$ at $\omega$ and its  minimum along $[2 \omega, 4 \omega]$ at $3\omega$.

\smallskip

\noindent Claim 6: $\Ima f(\omega) < \Ima f(\omega'+\omega) < \Ima f(\omega'+3\omega) < \Ima f(3\omega)$.

\smallskip

Claim 1 implies that $\Rea f$ is monotone along each of the named segments, and since $\Ima f$ differs between the two segments by Claim 2, 
$f$ must be one-to-one on the top side of $G$.
Claim 4 implies that $f$ is one-to-one on each of the two segments on the bottom side of $G$.
Claims 3, 5, and 6 imply that the images of these three segments do not intersect each other.
This shows that $f$ is one-to-one on the horizontal sides of $G$.

The claims can be established by the properties of $f'=F = \p_z v B$.
First note that, 
since $v(z)$ is positive in $G$ and vanishes on the horizontal sides of $G$, 
we have $\p_x v(z) = 0$ on both sides,
and for $x \in \RR$ we have $\p_y v(x+\omega_3) < 0$, and $\p_y v(x) > 0$.
In particular, $i \p_z v(z) = i (\p_x v - i \p_y v)/2 = \p_y v / 2$ is real.
The function $B(z)$ is a Jacobi $\text{sn}$ function, whose properties are well-known \cite[Section 47]{Akh}.
$B(z)$ sends the top side of $G$ to the unit circle, such that
the four segments $[\omega',\omega' + \omega]$, $[\omega' + \omega,\omega' + 2\omega]$, $[\omega' + 2\omega,\omega' + 3\omega]$, and $[\omega' + 3\omega,\omega' + 4\omega]$
correspond to the fourth, third, second, and first quadrants of the unit circle, respectively.
Multiplication by $\p_z v(z)$ distorts this circle and rotates it by an angle of $\pi/2$ 
(since $\p_z v(z)/i$ is positive) but preserves the two-fold symmetry.
This determines the sign of the real and imaginary parts of $f'$.
Since $dz=dx$ is purely real on the horizontal sides of $G$, this gives the monotonicity of $\Rea f$ stated in Claim 1.
Claims 2 and 3 follow from the sign of $\Ima f'$ and the fact that $\Ima f'$ is an odd function 
with respect to reflection in each of the points $\omega'+\omega$ and $\omega'+ 3\omega$.

The four segments $[0,\omega]$, $[\omega,2\omega]$, $[2\omega,3\omega]$, and $[3\omega,4\omega]$
on the bottom side of $G$ are sent to the second, first, fourth, and third quadrants of the unit circle respectively.
Since $\p_z v(z)/i$ is negative along the bottom side of $G$, under $f'(z)$ this becomes the first, fourth, third, and second quadrants, respectively.
This establishes Claim 4, and combined with the reflectional symmetry, also Claim 5.
Claim 6 follows from the fact that $\p_z v(z) B(z) > 0$ along the vertical segment $[\omega,\omega+\omega']$ and $\p_z v(z) B(z) < 0$ along $[3\omega,3\omega+\omega']$.

\smallskip

\noindent {\em Remark.} For the purpose of plotting Figure 1, instead of the above construction, 
we expressed $F$ as a ratio of Weierstrass $\sigma$ functions:
$$f'(z) = F(z) = \frac{\sigma(z-\omega+\omega'/2)^2 \cdot \sigma(z-3\omega+\omega'/2)^2 }{\sigma(z)^2 \cdot \sigma(z-2\omega) \cdot \sigma(z-6\omega+2\omega') }, $$
where $\sigma$ is a Weierstrass $\sigma$ function with fundamental ``periods'' $4\omega$, $2\omega'$ (but recall that $\sigma$ is not itself periodic).
As usual, 
the shifts are chosen based on the the zeros and poles of $F$,
but one of the shifts must be replaced by an equivalent lattice point in a different rectangle in order to 
satisfy \cite[Eq. (1), Sec. 14]{Akh}.
This explains why one of the poles is placed at $6\omega-2\omega'$.

\vspace{.1in}

\noindent
{\bf Example of type II.}
\vspace{.1in}

Only small modifications of the previous example are needed.
Using the same $G, G'$, $\omega_1,\omega_3$, we define $h$
as the $\omega_1$-periodic function, positive and harmonic in $G'$
except two logarithmic poles at $i\epsilon$ and $\omega_1+i\epsilon$,
where $\epsilon\in(0,\omega_3/2)$.
Then we can find constants $c_1$ and $c_2$ such that
$v=h+c_1y+c_2$ has critical points at $\omega_1/2+\omega_3/2$
and $3\omega_1/2+\omega_3/2$.

Then $v_z$ is an elliptic function with periods $\omega_1$, $2\omega_3$
with two simple poles at $i\epsilon$ and $-i\epsilon$ per period
parallelogram. This elliptic function has the form
$$\frac{-i\wp}{1+c\wp}+ic_0$$
with some small real $c$. The rest of the construction is the same
as in the previous example.

In a similar manner to the above, 
in order to plot the figures, 
we expressed $F$ as a ratio of Weierstrass $\sigma$ functions:
$$f'(z) = F(z) = \frac{\sigma(z-\omega+\omega'/2)^2 \cdot \sigma(z-3\omega +\omega'/2)^2 }{\sigma(z-i \epsilon) \cdot \sigma(z+i \epsilon) \cdot \sigma(z-2\omega-i\epsilon) \cdot \sigma(z-6\omega+i\epsilon +2\omega') }. $$

\begin{figure}[ht]
\centering
    \includegraphics[scale=.5]{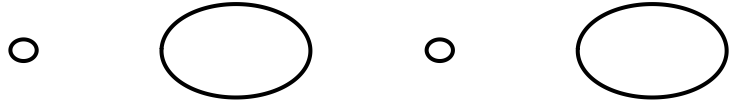}
    \caption{An example of Type II with $\omega_1=2$, $\omega_3=2$ and $\epsilon=0.5$.  
		Note that we have aligned the array horizontally in order to plot two periods.}
    \label{fig:2a}
\end{figure}

\begin{figure}[ht]
\centering
    \includegraphics[scale=.5]{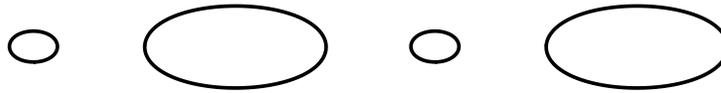}
    \caption{An example of Type II with $\omega_1=2$, $\omega_3=1.5$, and $\epsilon=0.4$.  
Note that we have aligned the array horizontally in order to plot two periods.}
    \label{fig:2b}
\end{figure}

\vspace{.1in}

\noindent
{\bf 6. Hollow vortex equilibria}
\vspace{.1in}

Let $G_j$ be smooth Jordan domains on the plane whose closures are
disjoint, and
$$D=\CC \setminus \cup_j D_j.$$
Let $F$ be a complex potential of a flow of an ideal fluid which is 
divergence-free and locally irrotational in $D$.
If the pressure (determined by $|F'|$ according to Bernoulli's law) 
is constant on $\partial D$ then $G_j$ can
be interpreted as constant-pressure gas bubbles in the flow.

The first examples of this situation, with two bubbles were constructed 
by Pocklington \cite{Pokl}. Periodic exceptional domains give
periodic examples with one bubble per period,
with the flow on the surface on the bubbles rotating in the same
direction \cite{1976} (see also \cite{Crowdy}). 
Crowdy and Green \cite{Crowdy}
constructed periodic examples with two bubbles per period rotating
in the opposite direction.
Our example of type II can be interpreted as a periodic flow with two
bubbles per period rotating in the same direction.

The velocity at infinity in our examples is directed in the opposite directions
on the two sides of the row of the bubbles.
\medskip

\noindent {\bf Acknowledgments:} We are grateful to Dmitry Khavinson for many helpful discussions
and to Razvan Teodorescu for a crucial observation regarding the construction of the examples of type I.
We also wish to thank Darren Crowdy for discussing with us the interesting connection to the hollow vortex problem.

{\em Department of Mathematics, Purdue University, West Lafayette, IN 47907 USA}
\end{document}